\documentclass{article}
\usepackage{latexsym}
\usepackage[]{amsfonts}
\usepackage{amsmath}

\title{\bf Strongly Embedded Subgroups of Groups of Odd Type}
\author{{\bf Christine Altseimer} \\ \\Institut f\"ur Mathematische
  Logik \\Albert-Ludwigs-Universit\"at
  Freiburg\\Eckerstr. 1\\ 79104 Freiburg\\ Germany}\date{} 

\begin{document} 
\pagestyle{plain}
\def\cc{\mathbb{C}}
\def\zz{\mathbb{Z}}
\def\ff{\mathbb{F}}
\def\nn{\mathbb{N}}
\def\qq{\mathbb{Q}}

\newtheorem{lemma}{Lemma}
\newtheorem{theorem}[lemma]{Theorem}
\newtheorem{corollary}[lemma]{Corollary} 
\newtheorem{proposition}[lemma]{Proposition}
\newtheorem{remark}[lemma]{Remark}
\newtheorem{definition}[lemma]{Definition}
\newtheorem{fact}[lemma]{Fact}
\newtheorem{con}{Conjecture}

\renewcommand{\refname}{\centerline{\bf References}}
\newcommand{\semi}{{\times}\kern-2pt \vrule height4.9pt width0.25pt depth0pt \kern2pt} 
\newcommand{\zs}{^{\circ}} 
\newcommand{\acf}{algebraically closed field}  
\newcommand{\fmr}{finite Morley rank} 
\newcommand{\uf}{if and only if}
\newcommand{\lemf}{\left(\begin{array}{cccc}} 
\newcommand{\rim}{\end{array}\right)}
\newcommand{\lemt}{\left(\begin{array}{cc}}

\def\H{\mathop{\rm H}\nolimits}
\def\G{\mathop{\rm G}\nolimits}
\def\K{\mathop{\rm K}\nolimits}
\def\PSL{\mathop{\rm PSL}\nolimits}
\def\PSp{\mathop{\rm PSp}\nolimits}
\def\SL{\mathop{\rm SL}\nolimits}
\def\GL{\mathop{\rm GL}\nolimits}

\def\proof{ {\it Proof.} $\,$}

\maketitle

\begin{abstract} \thispagestyle{empty}
 In this paper we prove that any strongly embedded subgroup of a
 $K^*$-group $G$ of \fmr\ and odd type that does not interpret any bad
 field and where $pr(G) \geq 2$ has to be solvable. If $n(G) \geq 3$ this has two important consequences. If there exists an involution $i \in G$ such that
 $C_G(i)$ is not solvable, then $G$ does not contain any proper
 2-generated core and centralisers of involutions have trivial cores.
\end{abstract}

\section{Introduction}
One of the great open problems in model theory is the Cherlin-Zil'ber
conjecture that states that an infinite simple group of \fmr\
is isomorphic, as an abstract group, to an algebraic group over an
v\acf. This paper belongs th the series of publications on the
classification of tame groups of \fmr. We call a group of \fmr\  {\it
  tame}\/ if none of its proper sections is a bad group and if it does
not interpret a bad field. Here a {\it bad group}\/ is a non-solvable group of \fmr\ all
of whose proper definable and connected subgroups are nilpotent and a
{\it bad field}\/ is a ranked structure of the form $\langle K,+,\cdot,
A\rangle$, where $K$ is an  \acf\ and $A$ is a proper infinite
multiplicative subgroup of $K^*$. 
Bad groups and bad fields are assumed not to exist, but the proof is
likely to require model theoretic methods. However it is hoped that
the classification of tame simple groups can be achieved using mainly ideas
from finite group theory. The first of these is to analyse a minimal
counterexample.  A group of \fmr\ is called a
$K$-{\it group}, if every infinite simple definable and connected
section of the group is an algebraic group over an algebraically
closed field.  A $K^*$-{\it group}\/ is a group of \fmr\ in which
every proper definable subgroup is a $K$-group.  In this inductive
setting the conjecture reduces to the following: 

\begin{con} \label{tame}
An infinite simple tame $K^*$-group is isomorphic to an algebraic
group over an algebraically closed field.  
\end{con}

 It is known \cite{mixed} that an infinite simple tame
$K^*$-group is either of {\it even type}, meaning that its Sylow
2-subgroups are nilpotent, definable and of bounded exponent, or of
{\it odd type}, meaning that the Sylow 2-subgroups are divisible
abelian by finite. In this paper we are dealing with the
classification of tame simple $K^*$-groups of odd type. For some
general background concerning the classification of tame simple groups
of odd type see \cite{nato}. We need to introduce some
terminology. If $E$ is a finite elementary 
abelian 2-group, its {\it 2-rank} $m(E)$ is the minimal number of
generators of $E$. If $H$ is any subgroup of a group $G$ of  \fmr\ and
odd type, its {\it 2-rank} $m(H)$ is the maximum of the 2-ranks of
elementary abelian subgroups in $H$. If $H$ is definable and $S$ is a
Sylow 2-subgroup of $H$, then the {\it normal 2-rank} $n(H)$ is the
maximum of the 2-ranks of normal elementary abelian subgroups in
$S$. 

Let $S$ be a Sylow 2-subgroup in a group $G$ of \fmr. We define the
{\it 2-generated core $\Gamma_{S,2}(G)$} as the definable closure of
the group generated by all normalizers $N_G(U)$ of all subgroup $U
\leq S$ with $m(U) \geq 2$.    

A group $G$ of \fmr\ will be called {\it quasi-simple}, if $G = G'$
and $G/Z(G)$ is non-abelian simple. By a {\it component} of $G$ we mean a connected definable
subnormal quasisimple subgroup. The {\it layer L(G)} of $G$ is the
product of all quasisimple subnormal subgroups of $G$. It is definable
and normal in $G$. The {\it Fitting subgroup} $F(G)$ is the group
generated by all normal nilpotent subgroups of $G$. $F(G)$ can be
shown to be nilpotent. The {\it generalized Fitting subgroup} $F^*(G)$
is taken to be $F\zs(G)*L\zs(G)$. Then $C\zs_G(F^*(G))\leq
F^*(G)$. $G$ satisfies the {\it B-conjecture}, if
$C\zs_G(i)=F^*(C_G(i))$ for any involution $i \in G$.  Simple
algebraic groups over an \acf\ of characteristic not 2 satisfy the
B-conjecture. Let $i \in G$ be an involution. A component $A \lhd
L(C_G(i))$ is called {\it intrinsic} if $i \in Z(A)$. An involution is
called {\it classical} if its centraliser contains an intrinsic
component isomorphic to $\SL(2,K)$ for an \acf\ $K$.

\begin{fact}
\label{three-disjunct}
Let $G$ be a simple tame $K^*$-group of odd type. Then one of the
following statements is true.
\begin{itemize}
\item $n(G) \leq 2$.
\item $G$ has a proper 2-generated core.
\item $G$ satisfies the $B$-conjecture and contains a classical involution.
\end{itemize}
\end{fact}

In this paper we show that the second case cannot occur, if $n(G) \geq
3$ and $G$ 
contains a non-solvable centraliser of an involution. For partial
results on the first case compare \cite{chr1} and \cite{chr2}. Berkman
\cite{berk} has furthermore done a nearly complete analysis of the
third case. 

Let $G$ be a group of \fmr. A proper definable subgroup $M$ of
$G$ is said to be {\it strongly embedded}\/ if it contains
 involutions and for every $g \in G\setminus M$, $M \cap M^g$ does not
 contain involutions. Alt{\i}nel has shown in \cite{alt94}, that any simple tame 
$K^*$-group of even type that has a strongly embedded subgroup is
isomorphic to $\SL_2(K)$ for an \acf\ $K$ of characteristic $2$. If
Conjecture \ref{tame} is true, then this is the only place, where
strongly embedded subgroups appear in the theory of simple tame
$K^*$-groups. The main aim of this section is to prove the following
partial result.

\begin{theorem}
\label{strong}
Let $G$ be a simple $K^*$-group of odd type that does not interpret any bad
field and such that $pr(G) \geq 2$. Let $M < G$ be a strongly embedded
subgroup. Then $M\zs$ is solvable. 
\end{theorem}

\begin{corollary}
Let $G$ be a simple tame $K^*$-group with a strongly embedded subgroup
$M$. Then $M$ is solvable or $G$ is of odd type and Pr\"ufer
2-rank $1$.
\end{corollary} 

\proof $G$ is either of odd or of even type by \cite{mixed}. If $G$ is of odd type and $pr(G)\geq 2$, then $M$ is solvable by Theorem
\ref{strong}. If $G$ is of odd type, then $M$ is solvable by
\cite{alt94}. \hfill $\Box$\vspace{1ex} 

Groups of \fmr\ with a strongly embedded
 subgroup have only one conjugacy class of involutions by
 \cite[10.19]{BN}.  Furthermore a proper subgroup $M$ of a group $G$ of
 \fmr\ is strongly embedded if and only if $M$ contains involutions,
 $C_G(t) \leq M$ for any $t \in I(M)$ and $N_G(S) \leq M$ for any
 Sylow 2-subgroup $S$ of $M$ by \cite{alt94}. Their structural
 properties are discussed in \cite{alt94} and slightly further
 explored in \cite{diplom}.

\section{Strongly embedded subgroups}

 We can describe strongly embedded
 subgroups differently, if $G$ is a connected $K^*$-group of \fmr\ and
 odd type such that $pr(G) \geq 2$.

\begin{proposition}
\label{embed}
Let $G$ be a connected $K^*$-group of \fmr\ and odd type such that
$pr(G) \geq 2$. Assume that $M < G$ is a strongly embedded subgroup of
$G$. Then 
\begin{itemize}
\item[(i)] $M\zs= \langle C_G(t)\zs |\; t \in D^* \rangle$ for any
four-subgroup $D$ of $M$  and
\item[(ii)]$M$ is a maximal proper definable subgroup of $G$. 
\end{itemize}
Thus, if\/ $G$ is a $K$-group, then it cannot contain a strongly
embedded subgroup.     
\end{proposition}

\proof Any Sylow 2-subgroup $S$ of $M$ is already a  Sylow
2-subgroups of $G$ by \cite{alt94}. Thus there
exists a four-subgroup $D \leq M$, as $pr(G) \geq 2$. Then $M\zs = \langle C_M(t)\zs |\;
t \in D^*\rangle$ by \cite[5.14]{nato}. On the other hand $C_G(t)\zs
\leq M\zs$ for all $t \in D^*$ as $M$ is strongly embedded, which gives us
that $M\zs = \langle C_G(t)\zs |\; t \in D^*\rangle$.  

We show that $M$ is a maximal proper definable subgroup of $G$. Let $N$ be a definable subgroup of $G$ such that $M \leq N < 
G$. Then $N$  is strongly embedded as well by
\cite{alt94}. Hence $N\zs = \langle C_G(t)\zs :
\; t \in D^*\rangle =M\zs$ by \cite[5.14]{BN} again. Especially $N \leq
N_G(M\zs)$. Let $g \in N_G(M\zs)$. Then $M \zs \leq M \cap M^g$
contains involutions as $I(M)
= I(M\zs)$ by \cite{alt94}. Since $M$ is strongly embedded, this
implies that $N_G(M\zs) \leq M$ and hence $N= N_G(M\zs)=M$.   

As finally $G = \langle C_G(t)\zs |\; t \in
D^*\rangle$ for any four-subgroup $D \leq G$, if $G$ is a $K$-group by
\cite[5.14]{nato}, $G$ cannot contain a strongly embedded subgroup in
this case. \hfill $\Box$   

\noindent\\ The reverse direction is true as well.

\begin{proposition}
\label{contra}
Let $G$ be a connected $K^*$-group of \fmr\ and odd type such that
$pr(G)\geq 2$. Let $N:=\langle C_G(l)\zs |\; l \in D^*
\rangle$ for a four-subgroup $D$ of $G$.  Then either $N$ is normal
in $G$ or $N_G(N)$ is a strongly embedded subgroup of $G$.
\end{proposition}

\proof Assume that $M:=N_G(N)<G$. Then $M$ is a $K$-group which
contains $D$ and $M\zs=N$ by \cite[5.14]{nato}. We are going to show
that $C_G(t) \leq M$ for any $t \in I(M)$ and that $N_G(S) \leq M$ for
any Sylow 2-subgroup $S$ of $M$. 

Let $t \in I(M)$ be any involution. Then $t$ is contained in Sylow
2-subgroup $S$ of $M$. As $Z(S) \neq 1$ by \cite[6.22]{BN}, there exists $s \in
S$ such that $D_1:=\langle t,s \rangle$  is a
four-subgroup of $S$ with $D_1 \cap Z(S) \neq 1$. Thus
$M\zs = \langle C_M(k)\zs|\; k \in D_1^* \rangle$ by
\cite[5.14]{nato}. Let $L:= \langle C_G(k)|\; k \in D_1^*
\rangle$. Then $C_G(t) \leq L$. We show that $L \leq M$. As $Z(S) \cap
D_1 \neq 1$, $S \leq L$. Let $g \in M$ such that $D^g \leq S$. Then 
\[L\zs =  \langle C_L(l)\zs|\; l \in (D^g)^* \rangle \leq
N^g=N.\] 
by \cite[5.14]{nato} again. Thus $M\zs \leq  L\zs \leq N = M\zs$,
$L\zs=N$ and $L \leq N_G(N) = M$. Hence $C_G(t) \leq M$ for any
involution $t \in M$.  

Let furthermore $S$ be any Sylow 2-subgroup of $M$. As $C_G(t) \leq M$
for all $t \in M$, $pr(M)=2$. As $N_G(S) \leq
N_G(S\zs)$, it is sufficient to show that $N_G(S\zs) \leq M$. Let $g
\in N_G(S\zs)$. Then either $g \in C_G(s) \leq M$ 
for some involution $s \in 
S\zs$ or there exists an elementary abelian subgroup $E \leq S\zs$
of order at least $4$ such that $g \in N_G(E)$. Since $M\zs =
\langle C_G(l)\zs:\; l \in E\backslash\{1\} \rangle$ by the first part
of the proof, this implies that $g \in N_G(M\zs)=M$.  
\hfill $\Box$

\section{Examples of strongly embedded subgroups}
We are now able to exhibit two possible examples of strongly
embedded subgroups.

\begin{remark}
\label{duppel}
Let $H$ be a $K$-group of odd type that contains an elementary abelian
2-subgroup $E$ of order $8$. Then 
\[H\zs=\langle C_H(D)\zs|\; D\leq E, [E:D]=2 \rangle.\]
\end{remark}

\proof   Let $D_1 \leq E$ be a four-subgroup. Then $H\zs=\langle
C_H(i) \zs | \; i \in D_1^*\rangle$ by \cite[5.14.]{nato}. However, $E
\leq C_H(i)$ for any $i \in D_1^*$ and thus 
\[C_H(i)\zs = \langle C_H(D)\zs | \; i \in D \leq E, [E:D]=2 \rangle\] 
by \cite[5.14.]{nato} again. \hfill $\Box$

\begin{proposition}\index{twog@$2$-generated core}
\label{twocore}
Let $G$ be a simple $K^*$-group of odd type such that $pr(G) \geq 2$
and assume that $G$ contains 
an elementary abelian subgroup $E$ of order $8$. If\/ $G$
contains a proper 2-generated core $\Gamma$, then $G$ contains a strongly
embedded subgroup $M:=N_G(\Gamma\zs)$. 
\end{proposition}

\proof We may assume that $\Gamma=\Gamma_{S,2}$, where $S$ is a Sylow
2-subgroup of $G$ that contains $E$. As $pr(G)\geq 2$, $S \leq
N_G(S\zs) \leq \Gamma$ by definition and $\Gamma\zs  = \langle
C_{\Gamma}(D)\zs|\; D\leq E, [E:D]=2 \rangle$ by Remark
\ref{duppel}. Since $N_G(D) \leq \Gamma$ for all four-subgroups $D
\leq E$,    
\[\Gamma\zs  = \langle C_G(D)\zs|\; D\leq E, [E:D]=2 \rangle =
\langle C_G(i) \zs|i \in D_1^* \rangle\] 
for any four-subgroup $D_1 \leq E$. The claim now follows by
Proposition \ref{contra}. \hfill $\Box$\vspace{1ex}

To give the second example, we need some more definitions. Let $G$ be a group of finite Morley rank. We say that $G$ is a
$2^{\bot}$-group if $G$ does not contain involutions and we denote by
$O(G)$ the {\it core of $G$}, the maximal normal definable connected
 $2^{\bot}$-subgroup of $G$.  The core exists since the product of two
 normal definable $2^{\bot}$-subgroups is a $2^{\bot}$-subgroup
 itself by \cite[ex.\ 11, p.\ 93]{BN}. Remember, that if $G$ is a
 K-group that does not interpret  a bad field, then $O(G)$ is
 nilpotent by \cite[5.8]{nato}. 

It is well known that in a simple algebraic group $G$ over
an \acf\ of characteristic 
not 2 we have $O(C_G(i))=1$ for any involution $i \in
G$. This should also hold in a simple tame $K^*$-group $G$
of odd type. In Corollary \ref{deadcores} we prove this conjecture in a
special case.

\noindent\\ Let $G$ be a group of \fmr\ and $\theta$ a function from
$I(G)$ into the set of definable subgroups of $G$.  We say that $\theta$
is a {\it signalizer functor for $G$}\/ if $\theta(s) \leq O(C_G(s))$ is a
connected definable normal subgroup of $C_G(s)$ for all $s \in I(G)$
and if furthermore for any commuting involutions $s,t \in G$ 
\[ \theta(t) \cap C_G(s) = \theta(s) \cap C_G(t).\]
In particular, $\theta$ is a signalizer functor for any simple $K^*$-group $G$
of odd type if $\theta(t) = O(C_G(t))$ for any involution $t \in G$ by
\cite[B.29]{BN}.  A signalizer
functor $\theta$ is called {\it complete}\index{signalizer functor!complete}, if for any elementary abelian
subgroup $E \leq G$ of order $\geq 8$ the subgroup $\theta(E) =
\langle \theta(t), \, t \in E \backslash \{1\} \rangle$ is a connected
$2^{\bot}$-subgroup and $C_{\theta(E)}(s) = \theta(s)$ for any $s \in
E^*$. We call $\theta$ {\it nilpotent}\/ if all the subgroups
$\theta(t)$ for $t \in I(H)$ are nilpotent.  

\begin{fact}[\cite{bor84}]
\label{signal}
Any nilpotent signalizer functor $\theta$ on a group $G$ of \fmr\ is
complete.
\end{fact}  

For the following proposition compare \cite[B.?]{BN} and \cite[9]{chr1}. 

\begin{proposition}
\label{trivcor}
Let $G$ be a simple $K^*$-group of odd type that does not interpret a
bad field. Assume that $pr(G) \geq 2$ and that $G$ contains an
elementary abelian subgroup $E$ of order $8$. Let $\theta$ be the
signalizer functor defined by $\theta(s)= O(C_G(s))$ for all $s \in
I(G)$. Then either $O(C_G(s))=1$ for all $s \in E^*$ or
$N_G\big(N_G(\theta(E))\big)$ is a strongly embedded subgroup of $G$.  
\end{proposition}

\proof As $G$ does not interpret a bad
field, $\theta$ is a nilpotent signalizer functor and thus complete by Fact
\ref{signal}. Hence $\theta(E)$ is nilpotent. Let $D$ be any four-subgroup of $E$. As $D$ acts on
$\theta(E)$ by a definable group automorphism  
\[ \theta(E) = \langle C_{\theta(E)}(t), \; t \in D^* \rangle \]
by \cite[4.6]{nato}. Since $\theta$ is complete this implies that 
\[ \theta(E) = \langle O(C_G(t)), \; t \in D^* \rangle \]            
In particular $N_G(D) \leq N_G(\theta(E))$. Thus   
\[N_G(\theta(E))\zs  = \langle C_G(D)\zs|\; D\leq E, [E:D]=2 \rangle =
\langle C_G(i) \zs|i \in D_1^* \rangle\] 
for any four-subgroup $D_1 \leq E$ by Remark \ref{duppel}. The claim
now follows by Proposition \ref{contra}. \hfill $\Box$

\section{Centralisers of involutions}
To prove Theorem \ref{strong}, we need the following crucial fact
about centralizers of involutions, which is also one of the main tools
to classify groups of small Pr\"ufer 2-rank.  

\begin{fact}[\cite{nato}]
\label{struccen}  
Let $G$ be a simple $K^*$-group of \fmr\ and odd type that does not
interpret a bad field. Let $i \in I(G)$. Then
$C_G(i)\zs/O(C_G(i))$ is a central product of an abelian
divisible group $T$ and a semisimple group
$H$ all of whose components are simple algebraic groups over algebraically
closed fields of characteristic different from 2. 
\end{fact}

\begin{lemma}
\label{basis} 
Let $H$ be a group of \fmr. For any $X \leq H$ and $h \in H$ let $\overline X:=
XO(H)/O(H)$ and $\overline t:=tO(H)$. 
\begin{itemize}
\item[(i)] If\/ $t \in I(H)$ then\/ $\overline t \in
  I(\overline{H})$ and if\/ $\overline s \in I(\overline{H})$ then
  there exists an involution in $sO(H)$. 
\item[(ii)] $O(\overline H)=1$.
\item[(iii)] Let $t, s \in I(H)$.  If\/ $[\overline t, \overline s]
  =1$, then there exists $u \in I(sO(H))$ such that
  $[t,u]=1$.
\item[(iv)] Let $t \in I(H)$.  Then $I(tO(H))=t^{O(H)}$. 
  Especially $t, s \in I(H)$ are conjugate in $H$ if and only if\/
  $\overline t, \overline s$ are conjugate in $\overline{H}$. 
\end{itemize}
\end{lemma}

\proof For the proof of $(i)$, $(iii)$ and $(iv)$ compare \cite[5]{chr1}.
To show $(ii)$, let $B \geq O(H)$ be such that $\overline B=O(\overline H)$. Then $B$ is
a definable, connected normal subgroup of $H$ since $O(H)$ and $O(\overline
H)$ are definable and connected. Furthermore $B$ cannot contain an
involution by $(i)$ and $B=O(H)$ by maximality of $O(H)$. \hfill $\Box$

\begin{lemma}
\label{sylowpruef}
Let $H$ be a group of \fmr\ and odd type. Let $S$ be a Sylow
2-subgroup of $H$ and $N \lhd H$ a definable solvable subgroup. Then
$S\zs N/N$ is the connected component of a Sylow 2-subgroup of $H/N$.
\end{lemma}

\proof Let $P$ be a Sylow 2-subgroup of $H/N$ that contains
$S\zs N/N$. Let furthermore $F$ be a subgroup of $H$ which
contains $N$ such that $F/N=d(P\zs)$. Then $F$ is a definable solvable
group and $S\zs \leq F$. Thus $S\zs$ is the connected component of a
Sylow 2-subgroup of $F$ and $S\zs N/N$ 
is the connected component of a Sylow 2-subgroup of $F/N=d(P\zs)$ by
\cite{hall}. Hence $S\zs N/N=P\zs$ and the claim follows.   \hfill $\Box$

\begin{lemma}
\label{prueferstruc}
Let $G$ be a simple $K^*$-group of \fmr\ and odd type that does not interpret
a bad field. Let $i \in G$ be an involution and $C:=C_G(i)\zs/O(C_G(i)) =H
*T$ as in Fact \ref{struccen}. If\/ $Q$ is a Sylow 2-subgroup of $H$ and
$S$ a Sylow 2-subgroup of $C$ that contains $Q$, then $S\zs = Q\zs R$ where
$R$ is the Sylow 2-subgroup of $T$. Especially  $pr(C)=pr(H)+pr(T)$,
were either $T$ is trivial or $pr(T)\geq 1$. 
\end{lemma}

\proof Since $C$ is the central product of $H$ and $T$, $H \cap T \leq
Z(H) \leq Z(C)$ and $H \cap T$ is finite, as $H$ is
semisimple. Let $Q$ be a Sylow 2-subgroup of $H$ and $S$ a Sylow
2-subgroup of $C$ that contains $Q$. Since $T$ is a central subgroup
of $C$, the Sylow 2-subgroup $R$ of $T$ is contained in $S$ and
$Q*R \leq S$. $T$ is as a divisible group connected by
\cite[ex.\ 3, p.\ 78]{BN}. Thus $R$ is connected as well by
\cite[9.29]{BN} and $Q\zs R \leq S\zs$. 

As $T$ is abelian $S\zs T/T$ is the connected component of a Sylow
2-subgroup of $C/T$ by Lemma \ref{sylowpruef} where $S\zs T/T \cong
S\zs/(S\zs \cap T) = S\zs/R$. On the other
hand  $Q\zs(H \cap T)/(H \cap T) \cong Q\zs/(Q\zs \cap T)$ is the connected
component of a Sylow 2-subgroup of $H/(H
\cap T)$  by Lemma \ref{sylowpruef}
again. As $C/T = HT/T \cong H/(H \cap T)$ actually $S\zs/R \cong Q\zs/(Q\zs
\cap T) = Q\zs/(Q\zs \cap R) \cong Q\zs R/R$ and  
$S\zs=Q\zs R$. As $pr(C)=pr(S)$, $pr(H)=pr(Q)$ and $pr(T) =pr(R)$ this implies
that $pr(C) =  pr(H) + pr(T)$.  

Assume now that $pr(T)=0$. As $R$ is connected, $R=1$. Especially $T
\leq O(C)$. However $O(C)=1$ by Lemma \ref{basis} and $T$ is
trivial. \hfill $\Box$

\begin{proposition}
\label{solvsyl}
Let $G$ be a simple $K^*$-group of odd type that does not interpret
a bad field. Then $N_G(P\zs)$ is solvable-by-finite for any Sylow
2-subgroup $P$ of $G$.
\end{proposition}

\proof Since $P\zs$ is a 2-torus $|N_G(P\zs)/C_G(P\zs)| <
\infty$ by \cite[6.16]{BN}. Hence it is enough to show that
$C_G(P\zs)\zs$ is solvable. 

Let  $i \in P\zs$ be any involution. As $P$ is infinite by
\cite{nato},  $i$ exists. We write $\overline X =
XO(C_G(i))/O(C_G(i))$ for any subgroup
$X \leq C_G(i)$.  If $C:=\overline{C_G(i)\zs}$ is solvable, the claim follows,
since $C_G(P\zs) \leq C_G(i)$. Assume now that $C$ is not
solvable. Then $C=H*T$ by Fact \ref{struccen}, where $T$ is an abelian
divisible group and $H$ is a nontrivial
semisimple group
all of whose components $H_n$ for $1 \leq n \leq k$  are simple
algebraic groups over algebraically
closed fields of characteristic different from $2$. $\overline P\zs$ is the
connected component of a Sylow 2-subgroup of $C$ by Lemma
\ref{sylowpruef}. Thus $\overline P\zs = Q\zs R$, where $Q$ is a Sylow
2-subgroup of $H$ and $R$ the Sylow 2-subgroup of
$T$ by Lemma \ref{prueferstruc}. Let $Q_n$ for $1 \leq n \leq k$ be
the projections of $Q\zs$ onto $H_n$. As $Q\zs$ consists of commuting
2-elements, $Q_n$
is an abelian 2-group and hence contained in a maximal torus $T_n$ of
$H_n$ for all $1 \leq n \leq k$ by \cite[15.4]{hum75}. As $Q_n$ are
the connected components of Sylow 2-subgroups of the simple algebraic
groups $H_n$ for all $1 \leq n \leq k$, $C_H(Q\zs) = T$, where
$T:=T_1\cdots T_k$. Thus $C_H(\overline P\zs) = T$ as well and
$\overline{C_{C_G(i)\zs}(P\zs)} \leq C_C(\overline P\zs)$ is abelian. Hence
$C_{C_G(i)\zs}(P\zs)\zs$ is a connected solvable group. However,
$C_G(P\zs)\zs \leq C_G(i)\zs$
and $C_G(P\zs)\zs \leq C_{C_G(i)\zs}(P\zs)\zs$ is solvable. \hfill $\Box$ 

\noindent\\ We need one more result to prove Theorem \ref{strong}.

\begin{lemma}
\label{solvcor}
Let $N$ be a group of \fmr\ and $H \lhd N$ a definable normal subgroup such that
$H\zs$ is a $2^{\bot}$-group. Then $N\zs=C_N(i)\zs H\zs$ for all
involutions $i \in H$. 
\end{lemma} 

Notice that the Lemma reduces to the Frattini argument
\cite[10.12]{BN}, if $\langle i \rangle$ is a Sylow 2-subgroup of
$H$. \vspace{1em} 

\proof  Let  $i \in I(H)$. Then
$H\zs=C_{H\zs}(i)H^-$ where $H^-$ is the
set of elements of $H\zs$ inverted by $i$ by \cite[ex.\ 14,
p.\ 73]{BN}. Furthermore 
\[(*) \quad rk(H\zs)=rk(C_{H\zs}(i)) + rk(H^-).\]  
We claim that $H^-=ii^{H\zs}=iI(iH\zs)$. Let $h \in H^-$. As $H\zs$ is
as a $2^{\bot}$-group 2-divisible by \cite[ex.\ 11, p.\ 72]{BN}, there
exists an element $d \in H\zs$ such that $h=d^2$. On the other hand
$h^i=h^{-1}$ and thus already $d=d^{-1}$ by \cite[ex.\ 12, p.\
72]{BN}. Thus $d \in H^-$ and $t=d^2=[i,d] \in ii^{H\zs}$. Furthermore
for all $d \in H\zs$, $ii^d \in iI(iH\zs)$ as $i^d=i[i,d] \in
iH\zs$. The claim follows since $iI(iH\zs)$ is a subset of $H\zs$
inverted by $i$.

$i^N$ on the other hand is the disjoint union of
finitely many sets $I(j_kH\zs)$ for $1 \leq k \leq d$ where $d
\leq |H/H\zs|$ and $j_k \in i^N$. Thus we may assume that
$rk(i^N)=rk(I(iH\zs))=rk(H^-)$ and 
\[(**) \quad rk(N) = rk(C_N(i)) + rk(i^N) = rk(C_N(i)) + rk(H^-). \] 
(*) and (**) imply that
\[rk(N) = rk(C_N(i))+ rk(H\zs)- rk(C_N(i) \cap H\zs) =
rk(C_N(i)H\zs).\]   
Hence $N\zs=C_N(i)\zs H\zs$ for any involution $i \in N$. \hfill
$\Box$\vspace{1ex}

\section{Proof of the theorem}

 Let $\sigma$ be the
solvable radical of $M\zs$. Since $M$ does not interpret a bad 
field $O(M\zs)$ is nilpotent and hence, as a definable connected
subgroup, contained in 
$\sigma\zs$. Notice that, since $O(M\zs)$ is characteristic in $M\zs$,
$O(M\zs) \leq O(M) \leq M\zs$ and thus $O(M\zs)=O(M)$. As furthermore
$M$ is a strongly embedded $K$-group
\[O(C_G(t))=O(C_M(t)) \leq O(M) \leq \sigma\zs \] 
for all involutions $t \in M$ by \cite{nato}. $\sigma\zs$ is a
characteristic subgroup of $M\zs$ and hence normal in $M$.  Since $M$
contains one conjugacy class of involutions as a strongly embedded
subgroup by \cite[10.19]{BN}, this implies that either $I(M) \subseteq
\sigma\zs$ or $\sigma\zs$ does not contain involutions. 

\noindent\\ {\bf $(a)$} {\it If\/ $\sigma\zs$ contains involutions,
  then $M\zs$ is solvable.} \vspace{1ex} 

\proof Assume that $I(M) \subseteq \sigma\zs$.  Let $S_1$ be a Sylow
2-subgroup of $\sigma\zs$ and
$S \geq S_1$ a Sylow 2-subgroup of $M$. Then $S$ is a Sylow 2-subgroup
of $G$ by \cite{alt94}. Furthermore $S_1$ is connected by
\cite[9.29]{BN} and thus a divisible subgroup of the abelian group
$S\zs$. Hence $S_1$ has a complement $B$ in $S\zs$ by the theorem of Baer \cite[15.1]{fuchs}.
As, however, all involutions are contained in $\sigma\zs$, $B$ 
has to be trivial, $S_1\zs =S\zs$ and $S\zs$ is a Sylow 2-subgroup of
$\sigma\zs$. By the Frattini argument
\cite[10.12]{BN} $M=N_M(S\zs)\sigma\zs$. Since 
$N_M(S\zs)/C_M(S\zs)$ is finite by \cite[6.16]{BN}, $M\zs =
C_M(S\zs)\zs \sigma\zs$. However $C_G(S\zs)\zs$ is solvable
by Corollary \ref{solvsyl}. Thus $M\zs$ is already solvable and $M\zs
= \sigma\zs$. \hfill $\Box$ 

\noindent\\ {\bf $(b)$} {\it If\/ $\sigma$ contains involutions, then
  either $I(M)$ is finite or $M\zs$ is solvable.} \vspace{1ex}  

\proof We may assume by {\bf $(a)$} that $\sigma\zs$ is a $2^{\bot}$-group and that $I(M) \subseteq \sigma$ as $\sigma$ is normal in
$M$. As $M\zs=C_G(i)\zs\sigma\zs$ for any involution $i \in M$ by Lemma
\ref{solvcor},  
\[\overline M\zs:=M\zs/\sigma\zs \cong C_G(i)\zs /(C_G(i)\zs \cap \sigma\zs) = C_G(i)\zs/O(C_G(i))\]   
for any involution $i \in M$.  Since $I(M)=I(M\zs)$ by \cite{alt94},  
$\{i\sigma\zs|i \in I(M)\} \subseteq Z(\overline M\zs)$. If $I(M)$ is
infinite, then $M\zs=\sigma$,
since $M\zs$ would otherwise contain infinitely many commuting
involutions by Lemma \ref{basis}. Hence $M\zs$ is solvable in this
case. \hfill $\Box$

\noindent\\ {\bf $(c)$} {\it If\/ $I(M)$ is finite, then either $pr(G) =1$ or
  $M\zs$ is solvable.}
\vspace{1ex}    

\proof If $I(M)$ is finite, then $M\zs = C_G(i)\zs$ for any involution
$i \in M$ as $rk(M) = rk(C_G(i)) + rk(I(M))$. Hence $M\zs/O(M)$ is a
central product of an abelian
divisible group $T$ and a semisimple group $H$ all of whose components are
algebraic groups over \acf s of characteristic not $2$ by
Fact \ref{struccen}. We assume that $M\zs$ is not solvable. Then
$O(M)=\sigma\zs$ by $(a)$ and $T=1$, as $T$ is connected.  Let $L_1,
\cdots, L_k$ be the components of $M\zs/O(M)=H$. By Lemma
\ref{basis} all involutions of $H$ are central in $H$. Especially
$|I(L_l)| = 2^{pr(L_l)-1}$, as $Z(L_l)$ is the intersection of all
maximal tori in $L_l$ for $1 \leq l \leq k$  by \cite[ex.\ 2,  p.\
162]{hum75} and $|I(H)| = 2^{pr(H)-1}$, as $Z(H)$ is the intersection
of all maximal tori in $H$ by \cite[ex.\ 2,  p.\ 162]{hum75} again. Let $g \in
M$. As $g$ acts as automorphism on $M\zs/O(M)$, for all $1 \leq l \leq
k$ there exists $1 \leq s \leq k$ such that $L_l^g= L_s$ by
\cite[7.1]{BN} and thus $I(L_l)^g = I(L_s)$. If $i_1 \in I(L_1)$, then
$I(M\zs/O(M))=i_1^M = \{I(L_1), \cdots ,I(L_k)\}$ since $M$, and hence
$M/O(M)$ by Lemma \ref{basis}, contains one conjugacy class of
involutions. Set $p:=pr(L_1)$. Then $pr(H)=kp$
and $|\{I(L_1), \cdots ,I(L_k)\}| \leq k 2^{p-1}$. This can only
happen for $k=1$, since
$|I(H)|=2^{kp-1} > k 2^{p-1} \geq |i_1^M|$ for $k>1$. Hence $H$ is
already quasi-simple. However, the only quasi-simple algebraic group in
which all involutions are central is $\SL_2(K)$ for an \acf\ of
characteristic not 2. (Compare e.g. \cite{lyons}). As
$pr(H)=pr(M\zs/O(M))=pr(M\zs)$ by Lemma \ref{sylowpruef}, the claim
follows. \hfill $\Box$  \vspace{1em}

Suppose that $\sigma$ does not contain involutions. Then $M\zs$ is
not solvable since $pr(M)=pr(G) \geq 2$ and $1 \neq \overline M:=M/\sigma$
contains one conjugacy class of involutions as in the proof of Lemma
\ref{basis}.  As 
furthermore $M$ is a K-group $\overline{M\zs}$ is
the direct sum of simple algebraic groups over algebraically closed
fields of characteristic not $2$ by \cite{alt94}. Let $L_1, \cdots,
L_k$ be the components of $\overline M\zs$. As $L_l$ contains involutions
for any $1 \leq l \leq k$, $\overline M$ acts transitively on the
components and all components are isomorphic by \cite[7.10]{BN}. 

We claim that $L_1$ contains one conjugacy class of
involutions. For assume that $i, j$ in $I(L_1)$ are not conjugate in
$L_1$. Then there exists $m \in \overline M$ such that $i^m = j$. As
$\overline{M\zs}$ is the direct product of conjugates of $L_1$, this implies
that $L_1^m=L_1$. Thus $C_{L_1}(i)^m=C_{L_1}(j)$. However,
non-conjugate involutions of simple algebraic groups have
non-isomorphic centralizers of involutions by \cite[4.3]{lyons}. Contradiction. 

The only simple algebraic groups with one conjugacy class of
involutions are $\PSL_2(K), \PSL_3(K)$ or $\G_2(K)$ for an \acf\ $K$
of characteristic not 2. Hence we may assume that $\overline
M\zs$ is isomorphic to a direct sum of copies of these groups.  On
the other hand $M\zs=C_G(i)\zs H$  for a connected $2^{\bot}$-group
$H$ which can be chosen independently from $i \in I(M)$ by
\cite[3.10]{alt94} since 
$M$ is strongly embedded. Hence $\overline M\zs$ is the product of the
centralizer of the involution $\overline{\mbox{\it \i}}$ and a Borel subgroup
$B$ of $\overline M\zs$ in the sense of algebraic group theory which
contains the connected nilpotent group $H\sigma/\sigma$.  As $H$ is independent
from $i\in I(M)$, we can choose $i$ such that $\overline{\mbox{\it \i}} \in
B$. Let $B=U \semi T$, 
where $U$ is the unipotent radical of $B$ and $T$ a maximal torus
which contains $\overline{\mbox{\it \i}}$ by
\cite[19.3]{hum75}. Let $U_l$ be
the projections of $U$ into $L_l$ and
$T_l$ be the projections of $T$ into $L_l$ for $1 \leq l \leq k$. Then
$\overline{\mbox{\it \i}} =i_1 \cdots i_k$ where $i_l \in I(T_l)$ and thus $T_l
\leq C_{L_l}(i_l)$  for $1 \leq l \leq k$. Hence $\dim L_l \leq
\dim C_{L_l}(i_l) + \dim U_l$ for $1 \leq l \leq k$. We may assume that
$i_1 \neq 1$.   

\noindent\\ If $L_1 \cong \PSL_2(K)$, then $\dim L_1  = 3$,
$\dim C_{L_1}(i_1) =1$ and $\dim U_1 =1$.

\noindent\\ If $L_1 \cong \PSL_3(K)$, then $\dim L_1  = 8$, $\dim C_{L_1}(i_1) =4$
and $\dim U_1 =3$.

\noindent\\ If  $L_1\cong \G_2(K)$, then $\dim L_1 = 14$,
  $\dim C_{L_1}(i_1) =6$ and $\dim U_1 =6$. 

\noindent\\ Thus none of these cases can occur. Contradiction. \hfill $\Box$

\section{Corollaries}

\begin{corollary} \index{twog@$2$-generated core}
\label{ftwocore}
Let $G$ be a simple $K^{*}$-group of odd type that does not interpret a
bad field. Assume that $pr(G) \geq 2$ and that $G$ contains an
elementary abelian subgroup of order $8$. If there exists an
involution $i \in G$ such that $C_G(i)$ is not solvable, then $G$ does
not contain a proper 2-generated core.
\end{corollary}

\proof Assume that $G$ does contain a proper 2-generated core. Then
$G$ contains a strongly embedded subgroup $M$ by Proposition
\ref{twocore}. If there exists an involution $i \in G$ such that
$C_G(i)$ is not solvable, we may assume that $i \in M$, as $G$
contains only one conjugacy class of involutions. Then $C_G(i) \leq M$
since $M$ is strongly embedded and $M$ is not solvable.  Contradiction
to Theorem \ref{strong}. \hfill $\Box$

\begin{corollary}\index{core!trivial}
\label{deadcores}
Let $G$ be a simple $K^{*}$-group of odd type that does not interpret a
bad field. Assume that $pr(G) \geq 2$ and that $G$ contains an
elementary abelian subgroup of order $8$. If there exists an
involution $i \in G$ such that $C_G(i)$ is not solvable, then
$O(C_G(t))=1$ for all $t \in E^*$. 
\end{corollary}

\proof By Proposition \ref{trivcor} either $O(C_G(t))=1$ for all $t
\in E^*$ or $G$ contains a strongly embedded subgroup $M$. The second
case cannot occur, however, as in the proof of Corollary
\ref{ftwocore}, which implies the claim. \hfill $\Box$   

\begin{corollary}
\label{genese}
Let $G$ be a simple $K^{*}$-group of odd type that does not interpret
a bad field. Assume that $pr(G) \geq 2$. If there exists an
involution $i \in G$ such that $C_G(i)$ is not solvable, then
$G=\langle C_G(k)\zs| \; k \in D_1^*\rangle$ for any four-subgroup
$D_1 \leq G$.   
\end{corollary}

\proof $D_1 \leq G$ be a four-subgroup of $G$ and set $M:=\langle
C_G(k)\zs| \; k \in D_1^*\rangle$. By Proposition \ref{contra}
either $G=N_G(M)$ or $M$ is a strongly embedded subgroup. Since
$G$ cannot contain a strongly embedded subgroup as in the proof of
Corollary \ref{ftwocore}, $G=M$ as $G$ is simple. \hfill $\Box$
\index{subgroup!strongly embedded|)}

\begin{center} \bf ACKNOWLEDGMENTS \end{center}
This paper is part of the author's Ph.D. thesis. She would like to
thank the Landesgraduiertenf\"orderung Baden-W\"urttemberg for their
financial support. She is very grateful to Alexandre V. Borovik, her
supervisor during her stay in Manchester, who gave all the help one
could wish for and without whom this work would not exist.

\end{document}